\documentclass[11pt]{article}
\usepackage{amsmath}
\usepackage{graphicx}
\usepackage[catalan,spanish,english]{babel} 
\usepackage[psamsfonts]{amssymb} 
\usepackage{cmmib57} 
\usepackage{exscale} 
\usepackage{amscd} 
\usepackage{latexsym} 
\usepackage{graphicx} 
\allowdisplaybreaks
\numberwithin{equation}{section}
\newtheorem{stat}{Statement}[section]
\newtheorem{theorem}[stat]{Theorem}

\newtheorem{lemma}[stat]{Lemma}

\newtheorem{remark}[stat]{Remark}
\newcommand{\IR}{{\mathbb R}}

\newcommand{\beq}{\begin{equation}}
\newcommand{\eeq}{\end{equation}}
\newcommand{\bal}{\begin{align}}
\newcommand{\eal}{\end{align}}
\newcommand{\beqn}{\begin{equation*}}
\newcommand{\eeqn}{\end{equation*}}
\newcommand{\baln}{\begin{align*}}
\newcommand{\ealn}{\end{align*}}
\newcommand{\tbar}{\bar t}
\newcommand{\xbar}{\bar x}

\begin{document}
%%%%%%%%%%%%%%%
\begin{titlepage}
%\hfill{\sl Draft version, February 7, 2008}
\vskip 1cm

\begin{center}
{\Large\bf Properties of the density 
for a three dimensional
stochastic wave equation}
\smallskip

by\\
\vspace{7mm}
{\sc Marta Sanz-Sol\'e}$\,^{(\ast)}$ \\
{\small Facultat de Matem\`atiques}\\
{\small Universitat de Barcelona} \\
{\small Gran Via de les Corts Catalanes 585}\\
{\small E-08007 Barcelona, Spain}\\
{\small marta.sanz@ub.edu}\\
{\small http://www.mat.ub.es/$\sim$sanz}
\end{center}
\vskip 2cm

\noindent{\bf Abstract.} We consider a stochastic wave equation in space dimension three driven by a noise white in time and with an absolutely continuous correlation measure given by the product of a smooth function and a Riesz kernel. Let $p_{t,x}(y)$ be the density of the law of the solution $u(t,x)$ of such an equation at points $(t,x)\in]0,T]\times \IR^3$. 
%We study the properties of the mapping $(t,x,y)\to p_{t,x}(y)$. In particular,
We prove that the mapping $(t,x)\mapsto p_{t,x}(y)$ owns the same regularity as the sample paths of the process $\{u(t,x), (t,x)\in]0,T]\times \IR^3\}$
established in  \cite{dss}. The proof relies on  Malliavin calculus and more explicitely,  the integration by parts formula of \cite{watanabe} and estimates derived form it.  
 \medskip

\smallskip

\noindent{\bf AMS Subject Classification}\\
{\sl Primary}: 60H15, 60H07.\\
{\sl Secondary}: 60G60, 60G17, 35L05.
\vskip 2cm

\noindent
\footnotesize
{\begin{itemize} \item[$^{(\ast)}$] Supported by the grant MTM 2006-01351 from the \textit{Direcci\'on General de
Investigaci\'on, Ministerio de Educaci\'on y Ciencia, Spain.}
\end{itemize}}
\end{titlepage}

\newpage

\section {\bf Introduction}

In this paper we consider the stochastic wave equation in space dimension $d=3$,
\begin{align} \label{0.1}
&\left(\frac{\partial^2}{\partial t^2} - \Delta \right) u(t,x) = \sigma\big(u(t,x)\big) \dot F(t,x) 
+ b \big(u(t,x)\big),\\ \nonumber
&u(0,x) = v_0(x),
\qquad \frac{\partial}{\partial t}u(0,x) = \tilde v_0(x),  
\end{align} 
where $t\in\, ]0,T]$ for a fixed $T>0$, $x\in\IR^3$, and 
$\Delta$ denotes the  Laplacian on $\IR^3$. The functions $\sigma$ and $b$ are Lipschitz continuous,
and the  process $\dot F$ is the {\it formal} derivative of a Gaussian random field, 
white in time and correlated in space defined as follows.

Let 
\beq
\label{0.2}
f(x)= \varphi(x)k_{\beta}(x),
\eeq
where $\varphi$ is a smooth positive function and $k_{\beta}$ denotes the Riesz kernel $k_{\beta}(x)= |x|^{-\beta}$, $\beta\in\, ]0,2[$.
We shall assume that $f$ defines a tempered measure and 
\beq
\label{cov}
\int_{\IR^3} \frac{\mu(d\xi)}{1+|\xi|^2} < \infty,
\eeq
where $\mu =\mathcal{F}^{-1}(f)$ and $\mathcal{F}$ denotes the Fourier transform operator. This condition is satisfied for instance for densities 
of the form (\ref{0.2}) with $\varphi(x)= \exp\left(-\sigma^2|x|^2/2\right)$
and $\beta\in]0,2[$ (see \cite{dss}).

Let $\mathcal{D}(\mathbb{R}^4)$ be the space of Schwartz test functions 
(see \cite{schwartz}).
Then, on some probability space, there exists a Gaussian process  $F = \left(F(\varphi),\ \varphi\in \mathcal{D}(\mathbb{R}^4)\right)$
 with mean zero and covariance functional defined by 
\begin{equation}\label{0.3}
E\big(F(\varphi) F(\psi)\big) = \int_{\mathbb{R}_+} ds \int_{\IR^3} dx  f(x)\left(\varphi(s)\ast\tilde\psi(s)\right)(x),
\end{equation}
where $\tilde\psi(s,x)= \psi(s,-x)$.

Riesz kernels are a  class of  singular correlation functions which have already appeared in several papers related with the stochastic heat and wave equations,
for instance in \cite{dalangfrangos}, \cite{dalang}, \cite{milletss}, \cite{milletss2}, \cite{milletmorien}. 
%They provide examples where condition (\ref{0.2}) is satisfied: for these covariances, (\ref{0.2}) is equivalent to the condition $0 < \beta < 2$ (see Example \ref{example2.5}).

%There are different possible approaches to giving a rigourous formulation to the initial value problem (\ref{0.1}). In all of them  
We recall that 
the fundamental solution $G(t)$ associated to the wave operator $\mathcal{L} = \frac{\partial^2}{\partial t^2} - \Delta$ in dimension three is given by 
%\label{0.4}
$G(t)= \frac{1}{4\pi t}\sigma_t,
$
$t>0$, where $\sigma_t$ denotes the uniform surface measure on the sphere of radius $t\in[0,T]$, hence with total mass $4\pi t^2$. 

The properties of $G(t)$ together with the particular form of the covariance of the noise play a crucial role in giving a rigourous formulation to the initial value problem (\ref{0.1}).
%%%%%%%%%%%

Here, we shall follow the same formulation as in \cite{dss} which 
for the purpose of existence and uniqueness of solution of (\ref{0.1})
introduces a localization of the SPDE by means of a set related with the past light cone,
as follows. Let ${\rm D}$ be a bounded domain in $\IR^3$. Set
\beq
\label{0.5}
K^{\rm D}_a(t) = \{y \in \IR^3: d(y,{\rm D})\leq a(T-t)\}, \ t \in [0,T],
\eeq
where $a\ge 1$ and $d$ denotes the Euclidean distance. Then, a solution to the SPDE (\ref{0.1}) in ${\rm D}$ is
an adapted, mean-square continuous stochastic process $\left(u(t) 1_{K^{\rm D}_a(t)},\ t \in [0,T] \right)$ with values in
$L^2(\IR^3)$, satisfying
\begin{align}
&u(t, \cdot) 1_{K^{\rm D}_a(t)}(\cdot) = 1_{K^{\rm D}_a(t)} (\cdot) \left(\frac{d}{dt} G(t) \ast v_0 + G(t) \ast \tilde v_0\right) (\cdot) \nonumber\\
&\qquad + 1_{K^{\rm D}_a(t)} (\cdot) \int_0^t \int_{\IR^3} G(t-s, \cdot -y) \sigma \left(u(s,y)\right) 1_{K^{\rm D}_a(s)}(y) M(ds, dy)
\nonumber\\
& \qquad + 1_{K^{\rm D}_a(t)} (\cdot) \int^t_0 ds\, G(t-s) \ast \left(b(u(s, \cdot))1_{K^{\rm D}_a(s)}(\cdot)\right) ,
\label{0.6}
\end{align}
a.s., for any $t \in [0, T]$, where we consider the stochastic integral defined in \cite{DM} and $M$ denotes the martingale measure
derived from $F$ (see \cite{dalangfrangos}).

The following result is a quotation of Theorem 4.11 in \cite{dss} and will be invoked repeatedly in this paper.
\begin{theorem}
\label{theorem1} Assume that:
\begin{itemize}
\item[(a)] the covariance density is of the form (\ref{0.2}) with the covariance factor $\varphi$ bounded and positive, $\varphi\in \mathcal{C}^1(\IR^d)$ and $\nabla \varphi \in \mathcal{C}_b^{\delta}(\IR^d)$, for some $\delta\in]0,1]$;
\item[(b)] the initial values $v_0$, $\tilde{v}_0$ are such that $v_0 \in\mathcal{C}^2(\IR^3)$, and $\Delta v_0$ and $\tilde v_0$ are H\"{o}lder continuous with orders $\gamma_1,\gamma_2 \in\, ]0,1]$, respectively;
\item[(c)] the coefficients $\sigma$ and $b$ are Lipschitz.
\end{itemize}
Then for any $q\in[2,\infty[$ and $\alpha \in\, ]0, \gamma_1 \wedge \gamma_2 \wedge \frac{2-\beta}{2}\wedge\frac{1+\delta}{2}[$,
there is $C >0$ such that for $(t,x), (\tbar,y) \in [0,T]\times {\rm D}$, 
\begin{equation}\label{0.7}
   E\left(\vert u(t,x) - u(\tbar,y)\vert^q\right) \leq C\left(\vert t - \tbar\vert + \vert x - y \vert\right)^{\alpha q}.
\end{equation}
In particular, a.s., the stochastic process $\left(u(t,x),\ (t,x) \in[0,T] \times {\rm D}\right)$ solution of (\ref{0.6}) has
$\alpha$-H\"{o}lder continuous sample paths, jointly in $(t,x)$, and
\beq
\label{bound1}
\sup_{(t,x)\in [0,T]\times {\rm D}} E\left(\vert u(t,x)\vert^q\right) < \infty
\eeq
for any $q\in[1,\infty[$.
\end{theorem}

In this paper we are interested in studying the properties of the density of the solution of (\ref{0.6}) as a function of 
$(t,x)\in]0,T]\times{\rm D}$, where ${\rm D}$ is a bounded subset of $\IR^3$.
 We shall denote this density by $p_{t,x}(y)$.  We shall prove that $(t,x)\to p_{t,x}(y)$ is jointly H\"older continuous, uniformly in $y$ on compact sets. 
% The Gaussian case

This question is trivial in the very particular case where the initial conditions $v_0$, $\tilde{v}_0$ and the coefficient
$b$ vanish, and the coefficient $\sigma$ is a constant function. In fact, with these assumptions and $\sigma=1$ the solution to the equation (\ref{0.6})
is a Gaussian process, centered, stationary in the space variable, and with
\begin{equation*}
\sigma_t^2:=E\vert u(t,x)\vert^2= \int_0^t ds \int_{\mathbb{R}^3} \mu(d\xi) \frac{\sin^2(s\vert\xi\vert)}{\vert\xi\vert^2}.
\end{equation*}
From the expression 
%\begin{equation*}
$p_{t,x}(y)= \frac{1}{\sqrt{2\pi}\sigma_t}\exp\left(-\frac{|y|^2}{2\sigma_t^2}\right)$
%\end{equation*}
it is not difficult to prove that
\begin{equation*}
\left\vert p_{t,x}(y)-p_{\bar{t},x}(y)\right\vert\le C (t_0,t_1, {\rm D})\left\vert t-\bar{t}\right\vert,
\end{equation*}
for $0<t_0\le t <\bar{t}\le t_1$, ${\rm D}\subset \mathbb{R}^3$.

However, in the  general situation that we are considering in this article, the problem becomes much more
involved.

Suppose that $v_0$, $\tilde{v}_0$ are null functions, assume also that the covariance of the 
process $F$ is given by (\ref{0.3}) with $dx f(x)$ replaced by $\Gamma (dx)$,
where $\Gamma$ is a non-negative, tempered, non-negative definite measure.
Set $\mu= \mathcal{F}^{-1}(\Gamma)$.
We introduce an assumption on $\mu$, denoted by (${\rm H}_\eta$), saying that  
$\int_{\IR^3} \frac{\mu(d\xi)}{(1+|\xi|^2)^\eta}<\infty$, for some value of $\eta\in]0,1]$. 
Assume that the coefficients $\sigma$ and $b$ are of class $\mathcal{C}^1$ with bounded derivatives and that 
(${\rm H}_\eta$) holds for some $\eta\in(0,1)$. Then, the existence of the  density  $p_{t,x}$ at any fixed point $(t,x)\in]0,T]\times {\rm D}$
has been established in \cite{qss1}.
Moreover, assuming that $\sigma$ and $b$ are $\mathcal{C}^\infty$ functions with bounded derivatives of order greater or equal to one, 
and that (${\rm H}_\eta$) holds for some $\eta\in(0,\frac{1}{2})$, it is proved in  \cite{qss2} that 
$y\mapsto p_{t,x}(y)$ is a $\mathcal{C}^\infty$ function. 
We refer the reader to \cite{ss} for results on applications of Malliavin calculus to the analysis of probability laws of
SPDEs.

In \cite{nq}, it is shown that the extension of Walsh's integral introduced in \cite{dalang} does not require for the integrands any
stationary property in the spatial variable. As a consequence of this fact, the results of \cite{dalang}, \cite{qss1}, \cite{qss2} and \cite{ss} concerning
the stochastic wave equation can be formulated with non null deterministic initial conditions. In addition, the solution of the equation in this setting
coincides with the solution to (\ref{0.6}).
Furthermore, in the particular case of absolutely continuous
covariance measures $\Gamma(dx)= f(x)dx$ satisfying (\ref{cov}) the existence and smoothness of the density $p_{t,x}$
are proved in \cite{nq} under the weaker assumption (${\rm H}_1$). 

Hence, on the basis of the above mentioned references and remarks, we can write
the next statement, which together with Theorem \ref{theorem1} are the starting point of our work.

\begin{theorem}
\label{t2} 
Assume assumptions (a) and (b) of Theorem \ref{theorem1}. Suppose also that  
$\sigma$ and $b$ are $\mathcal{C}^\infty$ functions with bounded derivatives of order greater or equal to one, and
$\inf\{|\sigma(z)|, z\in\IR\}\ge \sigma_0>0$.
Then, for any fixed $(t,x)\in]0,T]\times {\rm D}$, the law of the real valued random variable $u(t,x)$ solution to (\ref{0.6}) has a %$\mathcal{C}^\infty$
density $p_{t,x}\in \mathcal{C}^\infty$.
\end{theorem}

The main purpose of this paper is to prove that with the assumptions of this theorem, for any $y\in{\IR}$, the mapping
$$
(t,x)\in]0,T]\times \IR^3 \mapsto p_{t,x}(y)
$$
is jointly $\alpha$-H\"older continuous with $\alpha\in]0, \inf(\gamma_1,\gamma_2,\frac{2-\beta}{2},\frac{1+\delta}{2})[$ (see Theorem \ref{t1.1} in Section \ref{s1}).

For stochastic differential equations and some finite-dimensional stochastic evolution systems with an underlying semigroup structure one can find results of this type
for instance in (\cite{KS}). For SPDEs the problem has not been yet very much explored. To the best of our knowledge, this issue has only been studied for the stochastic heat equation in
spatial dimension $d=1$ in \cite{morien} and for the wave equation with $d=2$ in \cite{milletmorien} (see \cite{ballypardoux} and \cite{milletss} for  the existence and regularity of the density for these two types of SPDEs). 
It is worthy noticing that in these two references, the H\"older degree regularity of $p_{t,x}(y)$ in $(t,x)$ is better than for the sample paths of the solution process
$u(t,x)$, while in the equation under consideration we obtain the same order. As it will become clear from the proof,
the reason is the rather degenerate character of the fundamental solution of the wave equation in dimension three.

The method of our proof is based on the integration by parts formula of Malliavin calculus, as in the above mentioned references.
We next give the main ideas
and steps of the proof.
Fix $y\in\IR$ and let $(g_{n,y}, n\ge 1)$ be a sequence of smooth functions converging pointwise to the Dirac delta function $\delta_{\{y\}}$. Fix $t, \tbar \in]0,T]$, 
$x,\xbar\in {\rm D}$ and assume that we can prove 
\beq
\label{0.9}
\sup_{n\in\mathbb{N}, y\in K}\left|E\left[g_{n,y}\left(u(t,x)\right) - g_{n,y}\left(u(\tbar,\xbar)\right)\right]\right| \le C \left[|t-\tbar|^{\beta_1}+|x-\xbar|^{\beta_2}\right],
\eeq
for some $\beta_1, \beta_2>0$, where $K\subset \IR$. Then, since $p_{t,x}(y)= E\left[\delta_{\{y\}}\left(u(t,x)\right)\right]$ (see Theorem 1.12  in \cite{watanabe} for a rigorous meaning of this identity), by passing to the limit as $n\to\infty$, we will have joint H\"older continuity of the
mapping $(t,x)\in]0,T]\times {\rm D} \mapsto p_{t,x}(y)\in\IR$ with degree $\beta_1$ in $t$ and $\beta_2$ in $x$, uniformly in $y\in K$.

An estimate like (\ref{0.9}) is obtained by the following procedure. For simplicity we write $g$ instead of $g_{n,y}$. We first consider a Taylor expansion
of $g(u(\tbar,\xbar))$ around $u(t,x)$ up to a certain order $r_0$ chosen in such a way to obtain optimal values of $\beta_1$ and $\beta_2$. Then for any $r\le r_0$,
we estimate terms of the type
$$
\left| E\left[g^{(r)}\left(u(t,x)\right) \left(u(t,x)-u(\tbar,\xbar)\right)^r\right]\right|,
$$
and the term corresponding to the rest in the Taylor expansion, whose structure is similar. For this we use the version of the integration by parts formula
for one-dimensional random variables given in Lemma 2, page 54 of \cite{watanabe} (see also Equations (2.29)--(2.31) of \cite{nualart}) which we now quote as a lemma.
\begin{lemma}
\label{l1}
On an abstract Wiener space $(\Omega, \mathcal{\bf H}, P)$, we consider two real-valued random variables $\xi$ and $Z$ such that
$\xi\in \mathbb{D}^\infty$, $\Vert D\xi \Vert_{\mathcal{\bf H}}^{-1}\in \cap_{p\ge 2}L^p(\Omega)$, $Z\in \mathbb{D}^\infty$. Let $g$ be a function in $\mathcal{C}^r$,
for some $r\ge 1$. Denote by $\tilde{g}$ the antiderivative of $g$. Then,  the following formula holds:
\beq 
\label{0.10}
E\left(g^{(r)}(\xi) Z\right) = E\left(\tilde{g}(\xi)H_{r+1}(Z,\xi)\right),
\eeq
where $H_r$, $r\ge 1$, is defined recursively by
\begin{align*}
H_1(Z,\xi)&= \delta\left(\frac{Z D\xi}{\Vert D\xi\Vert^2_{\mathcal{H}}}\right),\\
H_{r+1}(Z,\xi)&= \delta\left(H_{r}(Z,\xi)\frac{Z D\xi}{\Vert D\xi\Vert^2_{\mathcal{H}}}\right), r\ge 1.
\end{align*}
\end{lemma}
In this Lemma, $\delta$ stands for the adjoint operator of the Malliavin derivative, also termed divergence operator
or Skorohod integral and we  have used the notations of \cite{nualart} and \cite{ss}, as we shall  do throughout the paper when referring to notions
and results of  Malliavin calculus.

The abstract Wiener space that we shall consider here is the one associated with the Gaussian process $F$ restricted to the time
interval $[0,T]$, as is described in  \cite{ss}, Section 6.1. For the sake of completeness and its further use, we recall that $\mathcal{\bf H}:=\mathcal{H}_T$,
$\mathcal{H}_T= L^2([0,T]; \mathcal{H})$ and that $\mathcal{H}$ is the completion of the inner product space consisting of 
test functions endowed with the inner product 
\beqn
\langle \varphi,\psi\rangle_{\mathcal{H}} = \int_{\IR^3} dx f(x) \left(\varphi\ast\tilde\psi\right)(x).
\eeqn

Assume that the function $\tilde{g}$ in Lemma \ref{l1} is bounded. From (\ref{0.10}) it clearly follows that
\beq
\label{0.11}
\left\vert E\left[g^{(r)}(\xi) Z\right]\right\vert \le \Vert \tilde{g}\Vert_\infty \Vert H_{r+1}(Z,\xi)\Vert_{L^1(\Omega)}.
\eeq
Furthermore, as a consequence of the continuity property of the Skorohod integral and the assumptions on $\xi$,
for any $r\ge 1$, $k\ge 1$ and $p\in(1,\infty)$,
$$\Vert H_r(Z,\xi)\Vert_{k,p} \le C\Vert Z\Vert_{k+r,4^r p}$$
(see Corollary 4.1 in \cite{morien}). Consequently, under the previous assumptions from (\ref{0.11}) we obtain
\beq
\label{0.12}
\left\vert E\left[g^{(r)}(\xi) Z\right]\right\vert \le C \Vert \tilde{g}\Vert_\infty \Vert Z\Vert_{r+1,4^{r+1}}.
\eeq
Let us recall that for a natural number $k$ and a real number $p\in[1,\infty[$,
\beqn
\Vert Z\Vert_{k,p}=\Vert Z\Vert_{L^p(\Omega)}+\sum_{r=1}^k\Vert D^{r}Z\Vert_{L^p\left(\Omega;\mathcal{H}_T^{\otimes r}\right)}.
\eeqn
We shall apply (\ref{0.12}) mainly to  $\xi:=u(t,x)$ and $Z:= \left (u(t,x)-u(\tbar,\xbar)\right)^r$, for natural values of $r$. 
Under the hypotheses of Theorem \ref{t1.1} the assumptions of Lemma \ref{l1} are satisfied (see \cite{qss2} and Chapters 7 and 8 of \cite{ss}). 
Thus we face the problem of giving upper bounds for $\Vert\left (u(t,x)-u(\tbar,\xbar)\right)^r\Vert_{r+1,4^{r+1}}$.

Malliavin derivatives of the solution of (\ref{0.1})  satisfy evolution equations (see \cite{ss}, Theorem 7.1 and \cite{qss2}). Indeed, for $x\in {\rm D}$, and a natural number $k\ge 1$, $(D_{.,*}^ku(t,x), (t,x)\in[0,T]\times {\rm D})$ is a $\mathcal{H}_T^{\otimes k}$-valued process satisfying $D_{\tau,*} u(t,x)=0$ for $\tau>t$,  and for $\tau\le t$ it is the solution of   the evolution equation
\begin{align}
&D_{\tau,*}^ku(t,x)= Z_{\tau,*}^k(t,x)\nonumber\\
&\quad+\int_0^t \int_{\IR^3}G(t-s,x-y)\nonumber\\
&\quad\quad\times \left(\Gamma^k(\sigma, u(s,y))+\sigma^\prime(u(s,y))D_{\tau,*}^ku(s,y)\right) M(ds,dy)\nonumber\\
&\quad+\int_0^tds \int_{\IR^3} G(s,dy) \left(\Gamma^k(b, u(s,y)) + b^\prime(u(t-s,x-y))D_{\tau,*}^ku(t-s,x-y)\right).\label{derivative}
\end{align}
In this equation, $\left(Z^k(t,x), (t,x)\in[0,T]\times D\right)$ is a $\mathcal{H}_T^{\otimes k}$--valued stochastic process and for a given function $g\in \mathcal{C}^k$
and a random variable $X\in\mathbb{D}^{k,2}$, $\Gamma^k(g,X)= D^k(g(X))-g^\prime(X)D^k X$. 

The solution of (\ref{derivative}) satisfies 
\beq
\label{norms}
\sup_{(t,x)\in[0,T]\times {\rm D}}\left\Vert u(s,y)\right\Vert_{k,p}<+\infty
\eeq
for any $p\in[1,\infty[$ (see Theorem 7.1 in \cite{ss}).

In the next section, we shall make use of the explicit form of
(\ref{derivative})
for $k=1$. In this case
\beq
\label{z1}
Z_{\tau,*}^1(t,x)= G(t-\tau,x-*)\sigma(u(\tau,*))
\eeq
and $\Gamma^1(g,X)=0$.

With some effort, using the tools on stochastic integration of Hilbert--valued processes developed in \cite{ss} it can be proved that the conclusions of Theorem \ref{theorem1} also apply to the 
$\mathcal{H}_T^{\otimes k}$--valued stochastic process solution to (\ref{derivative}). More precisely, for any $k\ge 1$, 
$q\in[2,\infty[$ and $\alpha\in]0,\gamma_1\wedge\gamma_2\wedge\frac{2-\beta}{2}\wedge\frac{1+\delta}{2}[$, there is $C>0$
such that for $(t,x), (\bar t,y)\in[0,T]\times D$,
\begin{equation}\label{hd}
\Vert u(t,x) - u(\tbar,y)\Vert_{k,q}\leq C\left(\vert t - \tbar\vert + \vert x - y \vert\right)^\alpha.
\end{equation}

Hence, with the H\"older continuity property on $u(t,x)$ and
its Malliavin derivatives we may be able to prove (\ref{0.9}) for specific values of $\beta_1, \beta_2$.

We shall fix what is the top order $r_0$ in the Taylor expansion of $g(u(\tbar.\xbar))$. Clearly, the  lower exponents $\beta_i$ should come from the first order term. However, in the examples studied so far, terms of first and second order give the same exponent. For the equation (\ref{0.6}) the situation is different. 
Already at the first order level of the expansion, we shall see that the contribution of the pathwise integral involving the coefficient $b$  is of the 
same order than
the H\"older continuity exponent given in Theorem \ref{theorem1}. Clearly, the second order term would provide twice the H\"older continuity degree. Therefore, a Taylor expansion of first order gives the best possible result. However, to conclude whether the regularity of the density $p_{t,x}$ in $(t,x)$ is the same as that of the sample paths of $u(t,x)$, we have to check that the contribution to the first order term in the Taylor expansion of the stochastic integral is not worse than that of the pathwise integral. This explains the strategy of
the proof of the main result in the next section.

%%%%%%%%%%%%%%%%%
%%%%%%%%%%%%%%%%%
%%%%%%%%%%%%%%%%% SECTION 1
\section{Main Result}
\label{s1}
Throughout this section ${\rm D}$  denotes a fixed bounded domain of $\IR^3$ and $C$ will be any positive finite constant. We assume that (\ref{cov}) holds. Our purpose is to prove the following theorem.

\begin{theorem}
\label{t1.1} Assume that:
\begin{itemize}
\item[(a)] the covariance density is of the form (\ref{0.2}) and the covariance density factor $\varphi$ is bounded and positive, $\varphi\in \mathcal{C}^1(\IR^d)$ and $\nabla \varphi \in \mathcal{C}_b^{\delta}(\IR^d)$
for some $\delta\in]0,1]$;
\item[(b)] the initial values $v_0$, $\tilde{v}_0$ are such that $v_0 \in C^2(\IR^3)$, and $\Delta v_0$ and $\tilde v_0$ are H\"{o}lder continuous with orders $\gamma_1,\gamma_2 \in\, ]0,1]$, respectively;
\item[(c)] the coefficients $\sigma$ and $b$ are $\mathcal{C}^\infty$ functions with bounded derivatives of order greater or equal to one,
 and there exist $\sigma_0>0$ such that $\inf\{|\sigma(z)|, z\in\IR\}\ge \sigma_0$.
\end{itemize}
Then the mapping
$$
(t,x)\in]0,T]\times {\rm D} \mapsto p_{t,x}(y)
$$
is $\alpha$-H\"older continuous jointly in $(t,x)$ with $\alpha\in]0, \inf(\gamma_1,\gamma_2,\frac{2-\beta}{2},\frac{1+\delta}{2})[$, uniformly in $y\in{\IR}$.
\end{theorem}

%%%%%%%%%%%
{\it Proof:} 
Fix $y\in \IR$ and let $(g_{n,y}, n\ge 1)$ be a sequence of regular functions converging pointwise to $\delta_{\{y\}}$ as $n\to\infty$; for example, a sequence of
Gaussian kernels with mean $y$ and variances converging to zero. We may assume that the corresponding antiderivatives $\tilde{g}_{n,y}$ are uniformly bounded
by $1$. To simplify the notation, we shall write $g$ instead of $g_{n,y}$.
\medskip

\noindent {\bf 1. Time increments}

\noindent  For $(t,x)\in[0,T]\times{\rm D}$ we consider the Taylor expansion 
\begin{align}
&E\left[g\left(u(t+h,x)\right)-g\left(u(t,x)\right)\right]=E\left[g^\prime\left(u(t,x)\right)\left(u(t+h,x)-u(t,x)\right)\right]\nonumber\\
&\quad +E\left[g^{\prime\prime}(\tilde u(t,x,h))\left(u(t+h,x)-u(t,x)\right)^2 \right],
\label{1.1}
\end{align}
where $h>0$ and $\tilde u(t,x,h)$ denotes a random variable lying on the segment determined by $u(t+h,x)$ and $u(t,x)$.
\smallskip

\noindent{\bf First order term}
\smallskip

\noindent Set 
\beqn
T_1(t,x,h) = \left\vert E\left[g^\prime\left(u(t,x)\right)\left(u(t+h,x)-u(t,x)\right)\right]\right\vert.
\eeqn
We aim to prove that
\beq
\label{theorem2}
\sup_{(t,x)\in[0,T]\times {\rm D}}T_1(t,x,h) \le C h^\alpha,
\eeq
with $\alpha\in ]0, \gamma_1 \wedge \gamma_2 \wedge \frac{2-\beta}{2}\wedge\frac{1+\delta}{2}[$.

By using equation (\ref{0.6}), we write $T_1(t,x,h)\le \sum_{i=1}^3 T_{1,i}(t,x,h)$, with
\begin{align*}
T_{1,1}(t,x,h)&=\Big|E\Big[g^\prime\left(u(t,x)\right)
\Big[\left(\frac{d}{dt} G(t+h) \ast v_0 + G(t+h) \ast \tilde v_0\right) (x)\\
& - \left(\frac{d}{dt} G(t) \ast v_0 + G(t) \ast \tilde v_0 \right)(x)\Big]\Big]\Big|\\
T_{1,2}(t,x,h)&=\Big|E\Big[g^\prime\left(u(t,x)\right)
\Big[\int_0^{t+h} ds \int_{\IR^3} G(t+h-s, dz) b(u(s,x-z))\\
& - \int_0^{t} ds \int_{\IR^3} G(t-s, dz) b(u(s,x-z))\Big]\Big]\Big|\\
T_{1,3}(t,x,h)&=\Big|E\Big[g^\prime\left(u(t,x)\right)\\
&\quad\times \Big[\int_0^{t+h} \int_{\IR^3} G(t+h-s, x-z) \sigma(u(s,x-z)) M(ds,dz)\\
& - \int_0^{t} \int_{\IR^3} G(t-s, x-z) \sigma(u(s,x-z))M(ds,dz)\Big]\Big]\Big|.
\end{align*}
In fact, by our choice of $(t,x)$ all the indicator functions in (\ref{0.6}) take the value $1$.

We shall  apply repeatedly the inequality (\ref{0.12}) with $r=1$, $\xi:= u(t,x)$
and different choices of $Z$. To begin with,
we take

$$Z:= \left(\frac{d}{dt} G(t+h) \ast v_0 + G(t+h) \ast \tilde v_0\right) (x)- \left(\frac{d}{dt} G(t) \ast v_0 + G(t) \ast \tilde v_0\right) (x).$$
Since $Z$ is deterministic, $\Vert Z\Vert_{k,p}=|Z|$, for any $k$ and $p$. Then, applying (\ref{0.12}) and Lemma 4.9 in \cite{dss} yields
\begin{equation}
\label{1.2}
\sup_{(t,x)\in]0,T]\times {\rm D}}T_{1,1}(t,x,h)\le C h^{\gamma_1\wedge \gamma_2}.
\end{equation}
We next study the term $T_{1,2}(t,x,h)$.
Let
$$
T_{1,2,1}(t,x,h)= \Big|E\Big[g^\prime\left(u(t,x)\right)\int_t^{t+h} ds \int_{\IR^3} G(t+h-s, dz) b(u(s,x-z))\Big]\Big|.
$$
We apply (\ref{0.12}) to $Z:= \int_t^{t+h} ds \int_{\IR^3} G(t+h-s, dz) b(u(s,x-z))$ and consider the measure on
$[t,t+h]\times \IR^3$ given by $ds G(t+h-s,dz)$ with total mass $\frac{h^2}{2}$ and an arbitrary $p\in[1,\infty[$. By applying Minkowski's inequality, we obtain
\begin{align*}
&\left\Vert \int_t^{t+h} ds \int_{\IR^3} G(t+h-s, dz) b(u(s,x-z))\right\Vert_{2,p}\\
&\le \int_t^{t+h} ds \int_{\IR^3} G(t+h-s, dz)\left\Vert b(u(s,x-z))\right\Vert_{2,p}
\end{align*}
By the chain rule of Malliavin calculus,
\beqn
\left\Vert b(u(s,y))\right\Vert_{2,p}\le C \Big(1+\left\Vert u(s,y)\right\Vert_{2,p}
+\left\Vert u(s,y)\right\Vert_{2,2p}^2\Big).
\eeqn
Consequently,
\beq
\label{b2}
\sup_{(s,y)\in[0,T]\times {\rm D}^{2T}}\left\Vert b(u(s,y))\right\Vert_{2,p} <\infty.
\eeq
where for a bounded set ${\rm D}\subset \mathbb{R}^3$ and $a\ge 0$, we denote by ${\rm D}^a=\{z\in \mathbb{R}^3; d(z,{\rm D})\le a\}$
and we have applied (\ref{norms}).

Thus, we have proved
\beqn
%\label{b2}
\sup_{(t,x)\in ]0,T]\times {\rm D}}\left\Vert \int_t^{t+h} ds \int_{\IR^3} G(t+h-s, dz) b(u(s,x-z))\right\Vert_{2,p}\le C h^2,
\eeqn
for any $p\in[1,\infty[$, and consequently, 

%\begin{align*}
%\left\Vert D^2 \left(b(u(s,x-z)\right)\right\Vert_{\mathcal{H}_T^{\otimes 2}}&\le \Vert b^{\prime\prime}\Vert_{\infty}
%\sup_{(s,z)\in[0,T]\times {\rm D}^{2T}}\Vert Du(s,z)\Vert^2_{\mathcal{H}_T}\\
%&+ \Vert b^{\prime}\Vert_{\infty}\sup_{(s,z)\in[0,T]\times {\rm D}^{2T}}\Vert D^2u(s,z)\Vert^2_{\mathcal{H}_T^{\otimes 2}},
%\end{align*}
\beq
\label{1.3}
\sup_{(t,x)\in ]0,T]\times {\rm D}}T_{1,2,1}(t,x,h)\le Ch^2.
\eeq
Set
$$
T_{1,2,2}(t,x,h)= \Big|E\Big[g^\prime\left(u(t,x)\right)\int_0^t ds \int_{\IR^3} b(u(s,x-z))[G(t+h-s,dz)-G(t-s,dz)]\Big]\Big|,
$$
that according to (\ref{0.12}) we can bound as follows,
$$
T_{1,2,2}(t,x,h)\le C \Big\Vert \int_0^t ds \int_{\IR^3} b(u(s,x-z))[G(t+h-s,dz)-G(t-s,dz)]\Big\Vert_{2,4^2}.
$$
We have 
\begin{align}
\label{1.3.bis}
&\int_0^t ds \int_{\IR^3} b(u(s,x-z))G(t+h-s,dz)=\int_0^t ds (t+h-s)\nonumber\\
&\quad\times\int_{B_1(0)} G(1,dz) b(u(s,x-(t+h-s)z)),\nonumber\\
&\int_0^t ds \int_{\IR^3} b(u(s,x-z))G(t-s,dz)=\int_0^t ds (t-s)\nonumber\\
&\quad\times\int_{B_1(0)} G(1,dz) b(u(s,x-(t-s)z)),
\end{align}
as can be easily checked by applying the 
change of variables $z\mapsto \frac{z}{t+h-s}$ and $z\mapsto \frac{z}{t-s}$, respectively.
Then, by the triangular inequality we obtain for any 
$p\in[1,\infty[$
\begin{align}
\label{b22}
%\sup_{(x,t)\in D\in]0,T]} 
&\left\Vert\int_0^t ds \int_{\IR^3} b(u(s,x-z))[G(t+h-s,dz)-G(t-s,dz)]\right\Vert_{2,p}\nonumber\\
& \le C h \left\Vert\int_0^t ds \int_{B_1(0)} G(1,dz) b\left(u(s,x-(t+h-s)z\right)\right\Vert_{2,p}\nonumber\\
&+ C\left\Vert\int_0^t ds \int_{B_1(0)} G(1,dz)\left[b\left(u(s,x-(t+h-s)z\right)-b\left(u(s,x-(t-s)z\right)\right]\right\Vert_{2,p}.
\end{align}
%and consequently,
%\beqn
%\sup_{x\in D}E\Big|\int_0^t ds \int_{\IR^3} b(u(s,x-z))[G(t+h-s,dz)-G(t-s,dz)]\Big]\Big|^p\le C \Vert G\Vert_\infty h^{\alpha p},
%\eeqn
%with $\alpha\in ]0, \gamma_1 \wedge \gamma_2 \wedge \frac{2-\beta}{2}\wedge\frac{1+\delta}{2}[$ and a constant independent of $t$ and $x$.
%%%%%%%First order derivative

For the study of the first term in the right-hand side of (\ref{b22}), we apply Minkowski's inequality and then (\ref{b2}). This yields
\beqn
\left\Vert\int_0^t ds \int_{B_1(0)} G(1,dz) b\left(u(s,x-(t+h-s)z\right)\right\Vert_{2,p}<C.
\eeqn
The Lipschitz property of $b$ and (\ref{0.7}), (\ref{hd}), (\ref{norms}) yield 
\beqn
\left\Vert b\left(u(s,x-(t+h-s)z\right)-b\left(u(s,x-(t-s)z\right)\right\Vert_{2,p}
\le C h^\alpha,
\eeqn
with $\alpha\in ]0, \gamma_1 \wedge \gamma_2 \wedge \frac{2-\beta}{2}\wedge\frac{1+\delta}{2}[$. 
Consequently, after having applied Minkowski's inequality we see that 
the second term of the right-hand side of (\ref{b22}) is bounded by $Ch^\alpha$, uniformly in $(t,x)\in[0,T]\times {\rm D}$.

Thus, we have proved 
\beqn
\sup_{(t,x)\in [0,T]\times {\rm D}}T_{1,2,2}(t,x,h) \le C  h^{\alpha},
\eeqn
and along with  (\ref{1.3}) we obtain
\beq
\label{1.4}
\sup_{(t,x)\in [0,T]\times{\rm D}}T_{1,2}(t,x,h) \le C  h^{\alpha},
\eeq
with $\alpha\in ]0, \gamma_1 \wedge \gamma_2 \wedge \frac{2-\beta}{2}\wedge\frac{1+\delta}{2}[$.

Let us remark that in \cite{morien} and  \cite{milletmorien} the contribution of the term analogous to $T_{1,2}(t,x,h)$
is a power of $h$ of higher order than the H\"older degree of the solution. In fact, for the heat equation and the wave equation in spatial dimension two,
 by integration of the increments $G(t+h-s,dz)-G(t-s,dz)$ we get powers of $h$. For the wave equation in dimension three,
such an approach  is not possible. Instead, ``increments" of $G(t-s,dz)$ are transfered to increments of $b(u(s,x-z))$
(this is the role played by the change of variables that we have performed to obtain (\ref{1.3.bis})) and after this, we can conclude by applying the Lipschitz property of $b$ 
and the H\"older continuity of the sample paths.

The inequality (\ref{1.4}) tell us that we are not going to improve the H\"older degree of the mapping
$t\in]0,T]\mapsto p_{t,x}(y)$ in more than the given $\alpha$. But it might happen that the contribution of $T_{1,3}(t,x,h)$
makes the overall estimation  worse. We next carry out  a careful analysis of this term and prove that its contribution
in terms of powers of $h$ is the same as $T_{1,2}(t,x,h)$.

We write $T_{1,3}(t,x,h)=T_{1,3,1}(t,x,h)+T_{1,3,2}(t,x,h)$ with
\begin{align*}
&T_{1,3,1}(t,x,h)=\Big|E\Big[g^\prime(u(t,x))\\
&\quad\times\int_t^{t+h}\int_{\IR^3} \sigma (u(s,y))G(t+h-s,x-y) M(ds,dy)\Big]\Big|,\\
&T_{1,3,2}(t,x,h)=\Big|E\Big[g^\prime(u(t,x))\\
&\quad\times\int_0^t \int_{\IR^3}\sigma (u(s,y))\left[G(t+h-s,x-y)-G(t-s,x-y)\right] M(ds,dy)\Big]\Big|.
\end{align*}
The term $T_{1,3,1}(t,x,h)$ vanishes, since the random variable $g^\prime(u(t,x))$ is adapted to the natural filtration
generated by the martingale measure $M$.%$\{M_t, t\ge 0\}$.

In contrast with $T_{1,2}$, for the analysis of $T_{1,3,2}(t,x,h)$ we do not start by applying (\ref{0.12}), which actually
would lead to worse results (see Remark \ref{r2.1}); instead,  
 we apply Proposition 3.9 of \cite{ss}. Since the mathematical expectation of a 
Skorohod integral is zero, we obtain
\begin{align*}
T_{1,3,2}(t,x,h)&=\Big|E\Big \langle g^{\prime\prime}(u(t,x)) D_{.,*}u(t,x), \sigma (u(.,*))\\
& \quad \times\left[G(t+h-.,x-*)-G(t-.,x-*)\right]1_{]0,t]}(\cdot) \Big\rangle_{\mathcal{H}_T}\Big|\\
&=\Big|E\Big[ g^{\prime\prime}(u(t,x))\Big \langle D_{.,*}u(t,x), \sigma (u(.,*))\\
& \quad \times\left[G(t+h-.,x-*)-G(t-.,x-*)\right]1_{]0,t]}(\cdot) \Big\rangle_{\mathcal{H}_T}\Big]\Big|.
\end{align*}
For any $(\tilde t,\tilde x)\in[0,T]\times \IR^3$, we define
\beqn
 B^h_{.,*}(\tilde t,\tilde x)= \sigma (u(.,*))\left[G(\tilde t+h-.,\tilde x-*)-G(\tilde t-.,\tilde x-*)\right],% 1_{]0,\tilde t]}(.),
 %Y^h(t,x)&=\left\langle D_{.,*}u(t,x),B^h(t,x)\right\langle.
 \eeqn
 With this notation, and by applying 
 (\ref{0.12}) to $T_{1,3,2}(t,x,h)$ we see that 
\beq
\label{t132}
T_{1,3,2}(t,x,h)\le C \left\Vert \left\langle D_{.,*}u(t,x),B^h{.,*}(t,x)\right\rangle_{\mathcal{H}_T}\right\Vert_{3,4^3},
\eeq
We shall consider norms $\Vert \cdot\Vert_{k,p}$ for arbitrary $k\ge 1$ and $p\in[1,\infty[$, instead of $\Vert \cdot\Vert_{3,4^3}$.
By virtue of (\ref{derivative}) and (\ref{z1}), we write
\begin{align*}
&\left\Vert\left\langle D_{.,*}u(t,x),B^h{.,*}(\tilde t,\tilde x)\right\rangle_{\mathcal{H}_T} \right\Vert_{k,p}\\
&\qquad \le  C \left(\Vert B_1^h(t,x;\tilde t,\tilde x)\Vert_{k,p} +\Vert B_2^h(t,x;\tilde t,\tilde x)\Vert_{k,p} +\Vert B_3^h(t,x;\tilde t,\tilde x)\Vert_{k,p}\right),
\end{align*}
where
\begin{align*}
B_1^h(t,x;\tilde t,\tilde x)&=\Big\langle G(t-.,x-*)\sigma (u(.,*)), B^h{.,*}(\tilde t,\tilde x)\Big\rangle_{\mathcal{H}_T},\\
B_2^h(t,x;\tilde t,\tilde x)&=\Big\langle \int_0^t \int_{\IR^3} G(t-s,x-y)\sigma^\prime(u(s,y)) D_{.,\ast} u(s,y) M(ds,dy),\\
&\quad\quad B^h{.,*}(\tilde t,\tilde x)\Big\rangle_{\mathcal{H}_T},\\
B_3^h(t,x;\tilde t,\tilde x)&=\Big\langle \int_0^t ds \int_{\IR^3} G(s,dy) b^\prime(u(t-s,x-y))D_{.,\ast} u(t-s,x-y), \\
&\quad\quad B^h{.,*}(\tilde t,\tilde x)\Big\rangle_{\mathcal{H}_T}.
\end{align*}
Consider the change of variables $(s,y)\to \left(t-s,\frac{y}{t-s}\right)$; 
Fubini's theorem along with Minkowski's inequality yield
\begin{align*}
&\left\Vert B_3^h(t,x;\tilde t,\tilde x)\right\Vert_{k,p}=\left\Vert  \int_0^t ds (t-s) \int_{\IR^3} G(1,dy) b^\prime(u(s,x-(t-s)y))\right.\\
&\quad\quad\times\left.\left\langle D_{.,*}u(s,x-(t-s)y),B^h{.,*}(\tilde t,\tilde x)\right\rangle_{\mathcal{H}_T} \right\Vert_{k,p}\\
&\quad\le C \Vert b^\prime\Vert_\infty \int_0^t ds \int_{\IR^3}  G(1,dy) \left\Vert\left\langle D_{.,*}u(s,x-(t-s)y),B^h{.,*}(\tilde t,\tilde x)\right\rangle_{\mathcal{H}_T} \right\Vert_{k,p}\\
&\quad\le C  \int_0^t ds \sup_{y\in K_a^D(s)}\left\Vert\left\langle D_{.,*}u(s,y),B^h{.,*}(\tilde t,\tilde x)\right\rangle_{\mathcal{H}_T} \right\Vert_{k,p} .
\end{align*}
The last inequality is obtained as follows. By definition of the sets $K_a^D(t)$, it is obvious that for any $(t,x)\in[0,T]\times D$, $x$ belongs to $K_a^D(t)$. Since the 
support of the measure $G(s,dy)$ is the boundary of the ball centered at zero and with radius $s$, the $y$--variable in the above integrals belongs to $K_a^D(s)$.

Hence, by Gronwall's lemma
\begin{align}
\label{g3}
&\sup_{x\in K_a^D(t)} \left\Vert\left\langle D_{.,\ast}u(t,x),  B^h_{.,\ast}(\tilde t,\tilde x)\right\rangle_{\mathcal{H}_T}\right\Vert_{k,p}\nonumber\\
&\qquad \le C \sup_{x\in K_a^D(t)} \left(\Vert B_1^h(t,x;\tilde t,\tilde x)\Vert_{k,p} +\Vert B_2^h(t,x;\tilde t,\tilde x)\Vert_{k,p} \right).
\end{align}

Our aim is to prove that
\beq
\label{g4}
\sup_{(t,x)\in[0,T]\times {\rm D}} \left\Vert\left\langle D_{.,*}u(t,x), B_{.,*}^h(t,x)\right\rangle_{\mathcal{H}_T}\right\Vert_{k,p}\le C h^\alpha,
\eeq
with $\alpha\in]0,\gamma_1\wedge \gamma_2\wedge\frac{2-\beta}{2}\wedge\frac{1+\delta}{2}[$, for  $k=3$. This is done recursively on 
$k=0,1,2,3$, where by convention $\Vert \cdot\Vert_{0,p}=\Vert \cdot\Vert_p$.

To illustrate the method and simplify the presentation, we shall consider in (\ref{g3}) the norm $\Vert\cdot\Vert_{1,p}$ instead of 
$\Vert\cdot\Vert_{3,p}$. That is, we shall deal only with derivatives up to the first order. Thus let as first prove (\ref{g4}) for $k=0$, that means for the $L^p(\Omega)$--norm.
For this, we start by studying
 the  $L^p(\Omega)$--norm of $B_2^h(t,x;\tilde t,\tilde x)$.
 
To shorten the notation, set $\mathcal{D}(s,y;\tilde t,\tilde x)=\langle D_{.,*}u(s,y),B^h_{.,*}(\tilde t,\tilde x)\rangle_{\mathcal{H}_T}$.
From Burkholder's inequality it follows that
\begin{align}
\label{g4.0}
&\Vert B_2^h(t,x;\tilde t,\tilde x)\Vert_p^p\nonumber\\
& =E\left|\int_0^t \int_{\mathbb{R}^3} G(t-s,x-y) \sigma^\prime\left(u(s,y)\right) \mathcal{D}(s,y;\tilde t,\tilde x) M(ds,dy)\right|^p\nonumber\\
&\le C E\Big(\int_0^t ds \int_{\mathbb{R}^3} dy \int_{\mathbb{R}^3} d\bar y G(t-s,x-y)\sigma^\prime\left(u(s,y)\right) f(y-\bar y)\nonumber\\
&\quad\times G(t-s,x-\bar y)\sigma^\prime\left(u(s,\bar y)\right)  \mathcal{D}(s,y;\tilde t,\tilde x)  \mathcal{D}(s,\bar y;\tilde t,\tilde x) \Big)^{\frac{p}{2}}\nonumber\\
&\le C \int_0^t ds \sup_{y\in K_a^D(s)}\left\{E\left\vert \left\langle D_{.,*}u(s,y),B^h_{.,*}(\tilde t,\tilde x)\right\rangle_{\mathcal{H}_T}\right\vert^p\right\}.
\end{align}
Therefore, (\ref{g3}) with $k=0$ and  Gronwall's lemma yields
\beqn
%\label{g5}
\sup_{x\in K_a^D(t)} \left\Vert\left\langle D_{.,\ast}u(t,x),  B^h_{.,\ast}(\tilde t,\tilde x)\right\rangle_{\mathcal{H}_T}\right\Vert_{p}
\le C\sup_{x\in K_a^D(t)} \left\Vert B_1^h(t,x;\tilde t,\tilde x)\right\Vert_{p}
\eeqn
with a constant $C$ independent of $\tilde t$ and $\tilde x$. Hence, we can fix
 $\tilde t=t$ and $\tilde x=x$ in the preceding inequality and obtain
 \beq
\label{g5}
\sup_{t\in[0,T]}\ \sup_{x\in K_a^D(t)} \left\Vert\left\langle D_{.,\ast}u(t,x),  B^h_{.,\ast}(t, x)\right\rangle_{\mathcal{H}_T}\right\Vert_{p}
\le C \sup_{t\in[0,T]}\  \sup_{x\in K_a^D(t)} \left\Vert B_1^h(t,x)\right\Vert_{p}
\eeq
where $B_1^h(t,x)$ stands for $B_1^h(t,x;t,x)$.

By the very definition of the inner product in $\mathcal{H}_T$ we have 
\begin{align*}
\left\Vert B_1^h(t,x)\right\Vert^p_p &= E\Big |\int_0^t dr \int_{\IR^3} d\xi \int_{\IR^3} d\eta G(t-r,x-\xi) \sigma(u(r,\xi)) f(\xi-\eta)\\
&\quad \times\big[G(t+h-r,x-\eta)
-G(t-r,x-\eta)\big] \sigma(u(r,\eta))\Big|^p.
\end{align*}
We consider each one of the terms in the difference of the right-hand side of this inequality and apply
respectively the change of variables
$$
 (\xi,\eta)\to\left(x-\xi,(x-\eta)\frac{t-r}{t+h-r}\right), \quad (\xi,\eta)\to(x-\xi,x-\eta).
$$
%that we apply to the first and second integrals of the right hand-side of the previous identity.
With this, the increments in time of the measure $G$ are transfered to increments of $\sigma$ and $f$.
More precisely, we obtain
$$
\left\Vert B_1^h(t,x)\right\Vert^p_p\le C(T_{1,3,2,1}(t,x,h)+T_{1,3,2,2}(t,x,h))
$$
with
\begin{align*}
&T_{1,3,2,1}(t,x,h)= E\Big |\int_0^t dr \int_{\IR^3} \int_{\IR^3}  G(t-r,d\xi)G(t-r,d\eta)\sigma (u(r,x-\xi))\\
&\quad \times\frac{t+h-r}{t-r} f\left(\frac{t+h-r}{t-r}\eta-\xi\right)
\left [\sigma(u(r,x-\frac{t+h-r}{t-r}\eta))-\sigma(u(r,x-\eta))\right]\Big|^p,\\
&T_{1,3,2,2}(t,x,h)= E\Big |\int_0^t dr \int_{\IR^3} \int_{\IR^3}  G(t-r,d\xi)G(t-r,d\eta)\sigma (u(r,x-\xi))\\
&\quad \times\sigma (u(r,x-\eta)) \left(\frac{t+h-r}{t-r}f\left(\frac{t+h-r}{t-r}\eta-\xi\right)-f(\eta-\xi)\right)\Big|^p .
\end{align*}
For the analysis of $T_{1,3,2,1}(t,x,h)$ we consider the measure with support on $[0,T]\times B_{t-r}(0)\times B_{t-r}(0)$ defined by
\beq
\label{n1}
\nu(dr;d\xi,d\eta):=dr\  G(t-r,d\xi) G(t-r,d\eta)
\left\vert\frac{t+h-r}{t-r}f\left(\frac{t+h-r}{t-r}\eta-\xi\right)\right\vert.
\eeq
Following the steps of the proof of Lemma 6.3 in \cite{dss} we obtain
$$
\sup_{0\le t\le t+h\le T}\int_0^t \int_{B_{t-r}(0)}  \int_{B_{t-r}(0)}  \nu(dr;d\xi,d\eta) <\infty.
$$
Then, we can write
\begin{align*}
T_{1,3,2,1}(t,x,h)&= 
E\Big\vert\int_0^t \int_{B_{t-r}(0)}  \int_{B_{t-r}(0)} 
\nu(dr;d\eta,d\xi) \sigma (u(r,x-\xi))\\
&\quad\times\left[\sigma(u(r,x-\frac{t+h-r}{t-r}\eta))-\sigma(u(r,x-\eta))\right]\Big\vert^p
\end{align*}
and apply H\"older's inequality with respect to the measure $\nu(dr;d\eta,d\xi)$. This yields
\begin{align*}
&T_{1,3,2,1}(t,x,h)\le C \int_0^t \int_{B_{t-r}(0)}  \int_{B_{t-r}(0)} \nu(dr;d\eta,d\xi)\\
&\quad\quad\times E\left\vert\sigma (u(r,x-\xi))\left[\sigma(u(r,x-\frac{t+h-r}{t-r}\eta))-\sigma(u(r,x-\eta))\right]\right\vert^p.
\end{align*}
We now apply Schwarz' inequality to the factor containing the expectation.
Since the coefficient $\sigma$ is a Lipschitz function, by (\ref{0.7}) and (\ref{bound1})  %and the support of the measure $G(t-r,dz)$ is the boundary of the ball $B_{t-r}(0)$,
we obtain
\begin{align}
\label{1.5}
T_{1,3,2,1}(t,x,h)&\le C\sup_{|\eta|=t-r}
 \left(E\left\vert\sigma(u(r,x-\frac{t+h-r}{t-r}\eta))-\sigma(u(r,x-\eta))\right\vert^{2p}\right)^{\frac{1}{2}}\nonumber\\
& \le C \sup_{|\eta|=t-r} \left\vert\frac{h}{t-r} \eta\right\vert^{\alpha p}
\le C h^{\alpha p},
\end{align}
with $\alpha\in ]0, \gamma_1 \wedge \gamma_2 \wedge \frac{2-\beta}{2}\wedge\frac{1+\delta}{2}[$ and a constant C not depending on $(t,x)\in[0,T]\times D$
for an arbitrary $D$.

To study $T_{1,3,2,2}(t,x,h)$ we consider the measure on $[0,t]\times B_{t-r}(0)\times B_{t-r}(0)$
given by
\begin{align*}
\mu^h(dr;d\xi,d\eta)&= dr\ G(t-r,d\xi)G(t-r,d\eta)\\
&\quad\times\left\vert \frac{t+h-r}{t-r}f\left(\frac{t+h-r}{t-r}\eta-\xi\right)-f(\eta-\xi)\right\vert.
\end{align*}
We also consider two additional measures with the same support as $\mu^h(dr;d\xi,d\eta)$ obtained by applying the triangular inequality to
the expression 
$$\left\vert \frac{t+h-r}{t-r}f\left(\frac{t+h-r}{t-r}\eta-\xi\right)-f(\eta-\xi)\right\vert.$$
 They are defined by
\begin{align}
\label{n2}
\mu_1^h(dr;d\xi,d\eta)&= dr G(t-r,\xi)G(t-r,d\eta)\frac{h}{t-r}  f(\eta-\xi),\nonumber\\
\mu_2^h(dr;d\xi,d\eta)&= dr G(t-r,\xi)G(t-r,d\eta)\frac{t+h-r}{t-r}\nonumber\\
&\quad\times\left\vert f\left(\frac{t+h-r}{t-r}\eta-\xi\right)-f(\eta-\xi)\right\vert.
\end{align}
With these new ingredients,
$
T_{1,3,2,2}(t,x,h)\le C\left(T_{1,3,2,2,1}+T_{1,3,2,2,2}\right),
$
where
\begin{align*}
T_{1,3,2,2,1}& = E\left|\int_0^t  \int_{B_{t-r}(0)} \int_{B_{t-r}(0)}\mu_1^h(dr;d\xi,d\eta) \sigma (u(r,x-\xi))\sigma (u(r,x-\eta))\right|^p,\\
T_{1,3,2,2,2}& = E\left|\int_0^t  \int_{B_{t-r}(0)} \int_{B_{t-r}(0)}\mu_2^h(dr;d\xi,d\eta) \sigma (u(r,x-\xi))\sigma (u(r,x-\eta))\right|^p.
\end{align*}

We next check that $\int_0^t\int_{B_{t-r}(0)} \int_{B_{t-r}(0)} \mu_1^h(dr;d\xi,d\eta) < Ch$. Indeed, owing to (\ref{0.2}) and by the change of
variable $r\mapsto t-r$, we have 
\begin{align*}
&\int_0^t\int_{B_{t-r}(0)} \int_{B_{t-r}(0)} \mu_1^h(dr;d\xi,d\eta)\\
&\le C h \int_0^t \frac{dr}{r} \int_{B_{t-r}(0)} \int_{B_{t-r}(0)} G(r,d\xi) G(r,d\eta)k_{\beta}(\xi-\eta)\\
&= C h \int_0^t \frac{dr}{r} \int_{\IR^3}  d\xi \frac{|\mathcal{F}G(r)(\xi)|^2}{|\xi|^{3-\beta}}
\le C h,
\end{align*}
uniformly in $t\in[0,T]$, where in the last inequality we have applied Lemma 2.3 of \cite{dss} with $b=1$.

Consequently, H\"older's inequality, the linear growth of the coefficient $\sigma$ and the property (\ref{bound1}) yield
\beq
\label{1.6}
T_{1,3,2,2,1}\le C h^p,
\eeq
%%%%%
%%%%%
%%%%% \mu_2
uniformly in $(t,x)\in[0,T]\times D$.

The next step consists of proving that $\mu_2^h(dr;d\xi,d\eta)$ defines a finite measure as well, and in giving an estimate of its total mass in terms of powers of $h$.
For this, we consider the inequality
\begin{align}
&\left\vert f\left(\frac{t+h-r}{t-r}\eta-\xi\right)-f(\eta-\xi)\right\vert \le \left\vert\varphi\left(\frac{t+h-r}{t-r}\eta-\xi\right)
-\varphi(\eta-\xi)\right\vert\nonumber\\
&\quad \times k_\beta(\eta-\xi)\nonumber\\
&\quad  +\left\vert\varphi\left(\frac{t+h-r}{t-r}\eta-\xi\right)\right| \left\vert k_\beta\left(\frac{t+h-r}{t-r}\eta-\xi\right)-k_\beta(\eta-\xi)\right\vert,\label{1.7}
\end{align}
which is a consequence of (\ref{0.2}) and the triangular inequality.
The properties of $\varphi$ together with Lemma 2.3 in \cite{dss} yield 
\begin{align}
&\int_0^t dr \int_{B_{t-r}(0)} \int_{B_{t-r}(0)} G(t-r,d\xi) G(t-r,d\eta) \frac{t+h-r}{t-r}\nonumber\\
&\quad\times \left\vert\varphi\left(\frac{t+h-r}{t-r}\eta-\xi\right)
-\varphi(\eta-\xi)\right\vert k_\beta(\eta-\xi)\nonumber\\
&\quad \le C h \int_0^t \frac{dr}{t-r} \int_{B_{t-r}(0)}
 \int_{B_{t-r}(0)} G(t-r,d\xi) G(t-r,d\eta) k_\beta(\eta-\xi)\nonumber\\
&\quad \le C h, \label{1.8}
\end{align}
uniformly in $t\in[0,T]$.

Let us now consider the contribution to $\mu_2^h(dr;d\xi,d\eta)$ of the second term of the right-hand side of (\ref{1.7}). Let  $\tilde\alpha\in]0,1[$, $\beta\in]0,2[$
with $\tilde\alpha+\beta\in]0,2[$. By applying Lemma 2.6 (a) of \cite{dss} with $b:=\tilde\alpha$, $a:=3-(\tilde\alpha+\beta)$, $c:=h$, $u:=\eta-\xi$, $x:=\frac{\eta}{t-r}$,
we obtain
\begin{align}
&\int_0^t dr \int_{B_{t-r}(0)} \int_{B_{t-r}(0)} G(t-r,d\xi) G(t-r,d\eta) \frac{t+h-r}{t-r}\left\vert\varphi\left(\frac{t+h-r}{t-r}\eta-\xi\right)\right\vert\nonumber\\
&\quad\times \left\vert k_\beta\left(\frac{t+h-r}{t-r}\eta-\xi\right)-k_\beta(\eta-\xi)\right\vert\nonumber\\
&\le\Vert\varphi\Vert_\infty \int_0^t dr \frac{t+h-r}{t-r} \int_{B_{t-r}(0)} \int_{B_{t-r}(0)} G(t-r,d\xi) G(t-r,d\eta) Dk_{\beta}\left(\eta-\xi,\frac{h}{t-r}\eta\right)\nonumber\\
&\le \Vert\varphi\Vert_\infty h^{\tilde\alpha} \int_0^t dr \frac{t+h-r}{t-r} \int_{B_{t-r}(0)} \int_{B_{t-r}(0)} G(t-r,d\xi) G(t-r,d\eta)\nonumber\\
&\quad \times \int_{\IR^3} dw k_{\tilde\alpha+\beta}(\eta-\xi-hw)\left\vert Dk_{3-\tilde\alpha}\left(w,\frac{\eta}{t-r}\right)\right\vert, \label{1.9}
\end{align}
where we have set $Dg(x,y):= g(x+y)-g(x)$ for a  function $g: \IR^3\to \IR$.

Our next purpose is to prove that the last integral in the above expression is bounded, uniformly in $t, t+h\in[0,T]$. For this, as in Lemma 6.4 of \cite{dss} we split the integral on the $w$-variable in the last expression into the sum of two
integrals: on a finite ball containing the origin and on the complementary of this set. In this way we obtain as an upper bound of 
\begin{align*}
&\int_0^t dr \frac{t+h-r}{t-r} \int_{B_{t-r}(0)} \int_{B_{t-r}(0)} G(t-r,d\xi) G(t-r,d\eta)\\
&\quad\times\int_{\IR^3} dw k_{\tilde\alpha+\beta}(\eta-\xi-hw)\left\vert Dk_{3-\tilde\alpha}\left(w,\frac{\eta}{t-r}\right)\right\vert,
\end{align*}
the sum of the three terms:
\begin{align*}
I_1^{(1)}& = \int_0^t dr \frac{t+h-r}{t-r} \int_{B_{t-r}(0)} \int_{B_{t-r}(0)} G(t-r,d\xi) G(t-r,d\eta)\\
&\quad\times \int_{B_2(0)} dw k_{\tilde\alpha+\beta}(\eta-\xi-hw)k_{3-\tilde\alpha}\left(w+\frac{\eta}{t-r}\right),\\
I_1^{(2)}& = \int_0^t dr \frac{t+h-r}{t-r} \int_{B_{t-r}(0)} \int_{B_{t-r}(0)} G(t-r,d\xi) G(t-r,d\eta)\\
&\quad\times\int_{B_2(0)} dw k_{\tilde\alpha+\beta}(\eta-\xi-hw)k_{3-\tilde\alpha}(w),\\
I_1^{(3)}& = \int_0^t dr \frac{t+h-r}{t-r} \int_{B_{t-r}(0)} \int_{B_{t-r}(0)} G(t-r,d\xi) G(t-r,d\eta)\\
&\quad\times \int_{B_2(0)^c} dw k_{\tilde\alpha+\beta}(\eta-\xi-hw)
\left\vert Dk_{3-\tilde\alpha}\left(w,\frac{\eta}{t-r}\right)\right\vert.
\end{align*}
%%%I_1^{(1)}

Consider the change of variable $w\mapsto w+\frac{\eta}{t-r}$ and then $\eta\mapsto \frac{t+h-r}{t-r}\eta$ that we apply to $I_1^{(1)}$.
By Fubini's theorem we obtain 
\begin{align*}
I_1^{(1)}& \le \int_{B_3(0)} dw\ k_{3-\tilde\alpha}(w) \int_0^t dr \frac{t+h-r}{t-r} \int_{B_{t-r}(0)} \int_{B_{t-r}(0)} G(t-r,d\xi) G(t-r,d\eta)\\
&\quad\times  k_{\tilde\alpha+\beta}\left(\frac{t+h-r}{t-r}\eta-\xi-hw\right)\\
&=\int_{B_3(0)} dw\ k_{3-\tilde\alpha}(w) \int_0^t dr \int_{\IR^3} \int_{\IR^3} G(t-r,d\xi) G(t-r+h,d\eta)\\
&\quad\times k_{\tilde\alpha+\beta}(\eta-\xi-hw).
\end{align*}

The properties of the Fourier transform and the expression of this operator applied to Riesz kernels yield, after regularization of $G$,
\begin{align*}
&\int_{\IR^3} \int_{\IR^3} G(t-r,d\xi) G(t-r+h,d\eta) k_{\tilde\alpha+\beta}(\eta-\xi-hw)\\
&=\int_{\IR^3}\mathcal{F}G(t-r)(\xi)\overline{\mathcal{F}G(t-r+h)(\xi)}k_{3-(\tilde\alpha+\beta)}(.-hw)(\xi).
\end{align*}
Hence by applying Schwarz's inequality, the last integral is bounded by
$$
\sup_{t\in[0,T]} \int_{\IR^3}|\mathcal{F}G(t)(\xi)|^2 k_{3-(\tilde\alpha+\beta)}d\xi,
$$
which is known to be finite whenever $\tilde\alpha+\beta\in[0,2]$ (see for instance Equation (2.5) in \cite{dss}).

Since $ k_{3-\tilde\alpha}(w)$ is integrable in a neighbourhood of the origin for any $\tilde\alpha>0$,
we finally obtain $I_1^{(1)}$, is bounded uniformly in $t,h \in[0,T].$

%%%I_1^{(2)}
Similar but simpler arguments show that the same property hold for $I_1^{(2)}$.

%%%I_1^{(3)}
For any $\lambda\in[0,1]$ set $\psi(\lambda)= k_{3-\tilde\alpha}\left(w+\lambda\frac{\eta}{t-r}\right)$. It is easy to check that
$|\psi^\prime(\lambda)|\le Ck_{4-\tilde\alpha}\left(w+\lambda \frac{\eta}{t-r}\right)$. Moreover, for $|w|\ge 2$, and $|\eta|= t-r$, $$\left\vert w+\lambda \frac{\eta}{t-r}\right\vert\ge \left\vert|w|-\left\vert \lambda \frac{\eta}{t-r}\right\vert\right\vert\ge |w|-1\ge \frac{|w|}{2},$$
by  the  triangular inequality. 

Thus,
\begin{align}
I_1^{(3)}& \le C \int_0^t dr \frac{t+h-r}{t-r}\int_{B_{t-r}(0)} \int_{B_{t-r}(0)} G(t-r,d\xi) G(t-r,d\eta)\nonumber\\
&\quad\times \int_{(B_2(0))^c} k_{\tilde\alpha+\beta}(\eta-\xi-hw) dw
\int_0^1 k_{4-\tilde\alpha}\left(w+\lambda \frac{\eta}{t-r}\right) d\lambda\nonumber\\
&\le C \left(\int_{(B_2(0))^c} dw k_{4-\tilde\alpha}(w)\right) \left(\int_0^t \frac{dr}{r}\int_{\IR^3}|\mathcal{F}G(r)(\xi)|^2 k_{3-(\tilde\alpha+\beta)} d\xi\right).\label{I}
\end{align}
For $\tilde\alpha\in]0,1[$ the integral $\int_{(B_2(0))^c}dw\  k_{4-\tilde\alpha}(w)$ is finite. Moreover, if $\tilde\alpha+\beta\in]0,2[$ the last integral in (\ref{I}) is also finite, owing to Lemma 2.3 in \cite{dss} applied to the
value $b=1$. This lead us to conclude that $I_1^{(3)}$ is bounded uniformly in $t,h \in[0,T].$

Summarizing, as a consequence of (\ref{1.8}), (\ref{1.9}) and the preceding discussion, we have proved that
\beq
\label{1.10}
\sup_{t\in[0,t]}\int_0^t\int_{B_{t-r}(0)}\int_{B_{t-r}(0)} \mu_2^h(dr;d\xi,d\eta)\le C h^{\tilde\alpha},
\eeq
with $\tilde\alpha\in ]0, \frac{2-\beta}{2}[$.

We can now apply H\"older's inequality with respect to the measure $\mu_2^h(dr;d\xi,d\eta)$. By virtue of (\ref{1.10}), the linear growth of $\sigma$ and (\ref{bound1}) we obtain
\beq
\label{1.11}
\sup_{(t,x)\in[0,T]\times {\rm D}}T_{1,3,2,2,2}\le C h^{\alpha p},
\eeq
with $\alpha\in ]0, \frac{2-\beta}{2}[$.

Finaly, the estimates (\ref{1.5}), (\ref{1.6}) and (\ref{1.11}) imply that 
\begin{equation}
\label{bp0}
\sup_{(t,x)\in[0,T]\times {\rm D}}\Vert B_1^h(t,x)\Vert_p\le C h^\alpha,
\end{equation}
and {\it a fortiori}
\begin{equation}
\label{bp}
\sup_{(t,x)\in[0,T]\times {\rm D}} \left\Vert\left\langle D_{.,*}u(t,x), B_{.,*}^h(t,x)\right\rangle_{\mathcal{H}_T}\right\Vert_p\le C h^\alpha,
\end{equation}
 with 
$\alpha\in ]0, \gamma_1 \wedge \gamma_2 \wedge \frac{2-\beta}{2}\wedge\frac{1+\delta}{2}[$. This finishes the analysis of the $\Vert\cdot\Vert_p$
contribution to the left-hand side of (\ref{g4}). 

%%%%%NORMA ||.||_{1,p}

We next consider the $L^p(\Omega, \mathcal{H}_T)$--norm of 
$D\left\langle D_{.,*}u(t,x), B_{.,*}^h(t,x)\right\rangle_{\mathcal{H}_T}$. As in the previous step, we shall replace $B_{.,*}^h(t,x)$ by $B_{.,*}^h(\tilde t,\tilde x)$ with arbitrary $\tilde t\in[0,T]$, $\tilde x\in\IR^3$.
By virtue of (\ref{g3}) and (\ref{bp}) it suffices to study the $L^p(\Omega, \mathcal{H}_T)$--norm of $DB_i^h(t,x;\tilde t,\tilde x)$ for $i=1,2$. We start with the analysis of $B_2^h(t,x;\tilde t,\tilde x)$.

By applying the differential rules of Malliavin calculus we have
\begin{align*}
&D_{.\ast} B_2^h(t,x;\tilde t,\tilde x) = G(t-.,x-\ast)\sigma^\prime(u(.\ast))\mathcal D(s,y;\tilde t,\tilde x)\\
&\quad+\int_0^t \int_{\IR^3}G(t-s,x-y)
\left[ \sigma^{''}(u(s,y))D_{.\ast}u(s,y)\mathcal D(s,y;\tilde t,\tilde x)\right.\\
&\quad \left.+\sigma^\prime(u(s,y))D_{.\ast} \mathcal D(s,y;\tilde t,\tilde x)\right]M(ds,dy).
\end{align*}
Applying H\"older's inequality and using that $\sigma^\prime$ is bounded, we obtain, as in (\ref{g4.0}),
\begin{align}
\label{g6}
&E\left\Vert G(t-.,x-\ast)\sigma^\prime(u(.\ast)) \mathcal D(s,y;\tilde t,\tilde x)\right\Vert_{\mathcal{H}_T}^p\nonumber\\ 
&\quad \le C \int_0^t ds \sup_{y\in K_a^D(s)}\left\{E\left\vert \mathcal D(s,y;\tilde t,\tilde x)\right\vert^p\right\}.
\end{align}

For fixed $\tilde t, \tilde x$ we consider the $\mathcal{H}_T$--valued process defined by
\beq
\label{g7}
K(s,y;\tilde t,\tilde x)= \sigma^{''}(u(s,y))D_{.\ast}u(s,y)\mathcal D(s,y;\tilde t,\tilde x)+\sigma^\prime(u(s,y))D_{.\ast} \mathcal D(s,y;\tilde t,\tilde x),
\eeq
$(s,y)\in[0,T]\times \IR^3$, for which we have 
\begin{align}
\label{g8}
E\left(\Vert K(s,y;\tilde t,\tilde x)\Vert_{\mathcal{H}_T}^p\right)&\le \Vert \sigma^{''}\Vert_{\infty}^p\left(E\left(\mathcal D(s,y;\tilde t,\tilde x)\right)^{2p}\right)^{\frac{1}{2}} 
\left(E\left\Vert Du(s,y)\right\Vert_{\mathcal{H}_T}^{2p}\right)^{\frac{1}{2}}\nonumber\\
&+\Vert \sigma^{\prime}\Vert_{\infty}^p E\left\Vert D\mathcal D(s,y;\tilde t,\tilde x)\right\Vert_{\mathcal{H}_T}^p.
%&\le C\left( h^{\alpha p} + E\left\Vert DB^h(s,y)\right\Vert_{\mathcal{H}_T}^p\right),
\end{align}
%where $\alpha\in]0, \gamma_1 \wedge \gamma_2 \wedge \frac{2-\beta}{2}\wedge\frac{1+\delta}{2}[$, and in the last inequality we have applied 
%(\ref{bp}).
We can apply the $L^p$--estimates for stochastic integrals with respect to the Gaussian process $M$ of Hilbert--valued integrands (see Equation (6.8) of
Theorem 6.1 in \cite{ss} and \cite{nq}, pg. 289) yielding
\begin{align}
\label{g9}
&\left\Vert \int_0^t \int_{\IR^3}G(t-s,x-y) K(s,y;\tilde t,\tilde x)
M(ds,dy)\right\Vert_{L^p\left(\Omega,\mathcal{H}_T\right)}^p\nonumber\\
&\qquad \le C\int_0^t ds \sup_{y\in K_a^D(s)} \left( \left(E\left(\mathcal D(s,y;\tilde t,\tilde x)\right)^{2p}\right)^{\frac{1}{2}} 
+\Vert D\mathcal D(s,y;\tilde t,\tilde x)\Vert_{L^p(\Omega;\mathcal{H}_T)}^p\right).
\end{align}
By taking $\tilde t=t$ and $\tilde x=x$ and considering the inequalities (\ref{g6}), (\ref{g9}), we obtain
\begin{align}
\label{g10}
&\Vert D B_2^h(t,x)\Vert^p_{L^p(\Omega;\mathcal{H}_T)}\le C\int_0^t ds \sup_{y\in K_a^D(s)}\left[E\left\vert\langle D_{.,\ast}u(s,y),B_{.,\ast}^h(t,x)\rangle_{\mathcal{H}_T}\right\vert^p\right.\nonumber\\
&\left.\quad + \left(E\left\vert\langle D_{.,\ast}u(s,y),B_{.,\ast}^h(t,x)\rangle_{\mathcal{H}_T}\right\vert^{2p}\right)^{\frac{1}{2}}
+\left\Vert D\langle D_{.,\ast}u(s,y),B_{.,\ast}^h(t,x)\rangle_{\mathcal{H}_T}\right\Vert_{L^p(\Omega;\mathcal{H}_T)}^p\right].
\end{align}
Then, (\ref{bp}) and  Gronwall's lemma yield
\begin{align}
\label{g11}
&\sup_{t\in[0,T]}\sup_{x\in K_a^{\rm D}(t)} \left\Vert D\left\langle D_{.,*}u(t,x), B_{.,*}^h(t,x)\right\rangle_{\mathcal{H}_T}\right\Vert_{L^p(\Omega;\mathcal{H}_T)}^p\nonumber\\
&\quad \le C \left(\sup_{t\in[0,T]}\sup_{x\in K_a^{\rm D}(t)} \left\Vert D B_1^h(t,x)\right\Vert_{L^p(\Omega;\mathcal{H}_T)}^p +h^{\alpha p}\right),
\end{align}
with $\alpha\in ]0, \gamma_1 \wedge \gamma_2 \wedge \frac{2-\beta}{2}\wedge\frac{1+\delta}{2}[$. 

The last step of the proof consist of checking that for an arbitrary bounded set ${\rm D}\subset \IR^3$,
\begin{equation}
\label{g12}
\sup_{(t,x)\in[0,T]\times {\rm D}} \left\Vert D B_1^h(t,x)\right\Vert_{L^p(\Omega;\mathcal{H}_T)}\le C h^{\alpha},
\end{equation}
with $\alpha\in ]0, \gamma_1 \wedge \gamma_2 \wedge \frac{2-\beta}{2}\wedge\frac{1+\delta}{2}[$. 

The proof of this fact can be done following the same lines as for (\ref{bp0}). We apply  the results on  the 
densities 
$\nu(dr;d\xi;d\eta)$, $\mu_1^h(dr;d\xi;d\eta)$, $\mu_2^h(dr;d\xi;d\eta)$, defined in (\ref{n1}), (\ref{n2}),
respectively, proved so far. Instead of the process $\{\sigma(u(s,y)), (s,y)\in[0,T]\times \IR^3\}$ and the
$L^p(\Omega)$--norm, we shall deal here with the $\mathcal{H}_T$--valued process $\{D\left(\sigma(u(s,y))\right), (s,y)\in[0,T]\times \IR^3\}$
and the $L^p(\Omega;\mathcal{H}_T )$--norm. In addition to (\ref{0.7}), we should also apply (\ref{hd}) and (\ref{norms}).
We leave the details to the reader.

%%%%%%%%%%%%%%%%%%%%%
%Hence, by (\ref{bp}), (\ref{g6}), (\ref{g9}) and Gronwall's lemma  we obtain
%\beq
%\label{g10}
%\sup_{(t,x)\in[0,T]\times {\rm D}}\Vert B^h(t,x)\Vert_{1,p}\le C h^\alpha.
%\eeq
%with $\alpha\in ]0, \gamma_1 \wedge \gamma_2 \wedge \frac{2-\beta}{2}\wedge\frac{1+\delta}{2}[$. 
Together with (\ref{1.2}) and (\ref{1.4}) this proves (\ref{theorem2}) and concludes the proof of the first step of the proof.
\medskip

\begin{remark}
\label{r2.1}
%The tricky approach to handle the term $T_{1,3,2}$ is justified by the following explanation.
Applying first (\ref{0.12}) and then estimates for the $\Vert \cdot\Vert_{2,p}$--norm of the stochastic integral
leads to
\begin{align*}
&T_{1,3,2}\\
&\quad \le C \left\Vert \int_0^t\int_{\IR^3} \sigma(u(s,y))\left[G(t+h-s,x-y)-G(t-s.x-y)\right]M(ds,dy)\right\Vert_{2,4^3}\\
&\quad \le C h^{\frac{\alpha}{2}}.
\end{align*}
Thus, we loose accuracy. This may be a justification of a pretty tricky approach in the preceding proof.
\end{remark}
\medskip

\noindent{\bf The rest in the time expansion}

The second and last term in (\ref{1.1}) to be examined is 
\beqn
R(t,x,h)= \left|E\left[g^{\prime\prime}(\tilde u(t,x,h))\left(u(t+h,x)-u(t,x)\right)^2 \right]\right|
\eeqn
We shall apply (\ref{0.12}) to the random variables $\xi:=\tilde u(t,x,h)$ and $Z:=\left(u(t+h,x)-u(t,x)\right)^2$. For this, 
we have to make sure that the assumptions of Lemma \ref{l1} are satisfied. For $Z:=\left(u(t+h,x)-u(t,x)\right)^2$,
and the two choices of $\xi$  -$u(t,x)$ and $u(t+h,x)$-  this has been proved in \cite{qss2}. Then it suffices to remark that
the norm $\Vert\cdot\Vert_{\mathcal{H}_T}$ as well as $\Vert\cdot\Vert_{\mathcal{H}_T}^{-1}$ define convex functions
and use the definition of $\tilde u(t,x,h)$ to conclude.

Consequently, 
\beqn
R(t,x,h)\le C \Vert u(t+h,x)-u(t,x)\Vert_{3,4^3}^2.
\eeqn
Owing to (\ref{0.7}) and (\ref{hd}) we conclude that
\beq
\label{1.12}
\sup_{(t,x)\in[0,T]\times {\rm D}}R(t,x,h)\le C h^{2\alpha}.
\eeq
with $\alpha\in ]0, \gamma_1 \wedge \gamma_2 \wedge \frac{2-\beta}{2}\wedge\frac{1+\delta}{2}[$.
\smallskip

The estimates (\ref{theorem2}) and (\ref{1.12}) show that
\beqn
\sup_{(t,x)\in[0,T]\times {\rm D}}\left\vert E\left[g(u(t+h,x))-g(u(t,x))\right]\right\vert \le C h^\alpha,
\eeqn
with $\alpha\in ]0, \gamma_1 \wedge \gamma_2 \wedge \frac{2-\beta}{2}\wedge\frac{1+\delta}{2}[$. Therefore
the mapping $t\in]0,T[\mapsto p_{t,x}(y)$ is H\"older continuous of degree
 $\alpha\in ]0, \gamma_1 \wedge \gamma_2 \wedge \frac{2-\beta}{2}\wedge\frac{1+\delta}{2}[$, uniformly in $y\in\IR^3$
 varying on bounded sets.
\medskip

\noindent{\bf Step 2: space increments}
\smallskip

Fix $t\in]0,T]$ and consider the Taylor expansion 
\begin{align}
&E\left[g\left(u(t,\xbar)\right)-g\left(u(t,x)\right)\right]=E\left[g^\prime\left(u(t,x)\right)\left(u(t,\xbar)-u(t,x)\right)\right]\nonumber\\
&\quad +E\left[g^{\prime\prime}(\hat u(t,x,\xbar))\left(u(t,\xbar)-u(t,x)\right)^2 \right],
\label{1.13}
\end{align}
where $x, \xbar \in D$ and  $\hat u(t,x,\xbar)$ denotes a random variable lying on the segment determined by $u(t,\xbar)$ and $u(t,x)$.
\smallskip

\noindent{\bf First order term}
\smallskip

Our aim is to prove that
\beq
\label{s}
\sup_{t\in[0,T]}\left\vert E\left[g^\prime\left(u(t,x)\right)\left(u(t,\xbar)-u(t,x)\right)\right]\right\vert \le C |x-\xbar|^\alpha,
\eeq
with  $\alpha\in ]0, \gamma_1 \wedge \gamma_2 \wedge \frac{2-\beta}{2}\wedge\frac{1+\delta}{2}[$.

As for the time increments, we consider Equation (\ref{0.6}) and write 
$$E\left[g^\prime\left(u(t,x)\right)\left(u(t,\xbar)-u(t,x)\right)\right]=\sum_{i=1}^3 S_{i}(t,x,\xbar),$$
 with 
\begin{align*}
S_{1}(t,x,\xbar)&=\Big|E\Big[g^\prime\left(u(t,x)\right)
\Big[\left(\frac{d}{dt} G(t) \ast v_0 + G(t) \ast \tilde v_0\right) (\xbar)\\
& - \left(\frac{d}{dt} G(t) \ast v_0 + G(t) \ast \tilde v_0 \right)(x)\Big]\Big]\Big|\\
S_{2}(t,x,\xbar)&=\Big|E\Big[g^\prime\left(u(t,x)\right)\\
&\quad\times
\int_0^{t} ds \int_{\IR^3} G(t-s, dz) \big[b(u(s,\xbar-z))-b(u(s,x-z))\big]\Big]\Big|\\
S_{3}(t,x,\xbar)&=\Big|E\Big[g^\prime\left(u(t,x)\right)\\
&\quad\times \int_0^{t} \int_{\IR^3} \big[G(t-s, \xbar-z)-G(t-s,x-z)\big]
\sigma(u(s,y)) M(ds,dz)\Big]\Big|.
\end{align*}

Let us consider $S_{1}(t,x,\xbar)$. As for the term $T_{1,1}(t,x,h)$, we first apply the inequality (\ref{0.12})
and notice that 
\beqn
Z(t;x,\bar x):= \left(\frac{d}{dt} G(t) \ast v_0 + G(t) \ast \tilde v_0\right) (\xbar)
 - \left(\frac{d}{dt} G(t) \ast v_0 + G(t) \ast \tilde v_0 \right)(x)
\eeqn
is deterministic. Thus, it suffices to estimate the absolute value of the random variable $Z(t;x,\bar x)$ defined before.
For this, we apply Lemmas 4.2 and 4.4 in \cite{dss} which tell us that the fractional Sobolev norm of any integration degree
$p\ge 2$  and differential order $\rho< \gamma_1\wedge\gamma_2$ is bounded. Hence, since $p$ is arbitrary, 
by the Sobolev embedding theorem
we have that
\beq
\label{1.14}
\sup_{t\in[0,T]}S_{1}(t,x,\xbar)\le \sup_{t\in[0,T]}C |Z(t;x,\bar x)|\le C |x-\xbar|^\rho,
\eeq
with $\rho<\gamma_1\wedge\gamma_2$.

%%%%%%%%%%%%%%
%%%%%%%S_2(t,x,\bar x)
%%%%%%%%%%%%%%
We continue the proof with the study of the term $S_{2}(t,x,\xbar)$. By virtue of (\ref{0.12}), it suffices to find an upper bound
of 
\beqn
\left\Vert \int_0^{t} ds \int_{\IR^3} G(t-s, dz) \big[b(u(s,\xbar-z))-b(u(s,x-z))\big]\right\Vert_{2,4^2}
\eeqn
in terms of a power of $|x-\xbar|$.  

The measure on $[0,t]\times \IR^3$ defined by $ds\ G(t-s,dz)$ is finite. Hence, we can apply Minkowski's  inequality
and obtain for any $p\in[1,\infty[$
\begin{align*}
&\left\Vert\int_0^t ds \int_{\IR^3} G(t-s, dz) \left[b(u(s,\xbar-z))-b(u(s,x-z))\right]\right\Vert_{2,p}\\
&\quad \le \int_0^{t} ds \int_{\IR^3} G(t-s, dz) \left\Vert b(u(s,\xbar-z))-b(u(s,x-z))\right\Vert_{2,p}\\
&\quad \le C |x-\xbar|^\alpha,
\end{align*}
with $\alpha\in ]0, \gamma_1 \wedge \gamma_2 \wedge \frac{2-\beta}{2}\wedge\frac{1+\delta}{2}[$. The last inequality is obtained by using that  
 $b$ and its derivatives are Lipschitz continuous and bounded functions, and by applying (\ref{0.7}) and (\ref{hd}).

Hence,
\beq
\label{1.15}
\sup_{t\in[0,T]} S_{2}(t,x,\xbar)\le C |x-\xbar|^\alpha,
\eeq
for any $\alpha\in ]0, \gamma_1 \wedge \gamma_2 \wedge \frac{2-\beta}{2}\wedge\frac{1+\delta}{2}[$.

%%%%%%%%%%%%%%
%%%%%%%S_3(t,x,\bar x)
%%%%%%%%%%%%%%

To analyze $S_{3}(t,x,\xbar)$ we proceed in a similar manner as for $T_{1,3}(t,x,h)$ by applying first
Proposition 3.9 in \cite{ss} and then (\ref{0.12}).  We obtain
\begin{align*}
&S_{3}(t,x,\xbar)=\Big|E\Big[g^{\prime\prime}(u(t,x))\langle D_{.,\ast}u(t,x), \sigma(u(.,\ast))\\
&\quad\times\left[G(t-.,\xbar-\ast)-G(t-.,x-\ast)\right]1_{]0,t]}(\cdot)\rangle_{\mathcal{H}_T}\Big]\Big|\\
&\le C \Big\Vert \langle D_{.,\ast}u(t,x), \sigma(u(.,\ast))
\left[G(t-.,\xbar-\ast)-G(t-.,x-\ast)\right]1_{]0,t]}(\cdot)\rangle_{\mathcal{H}_T}\Big\Vert_{3,4^3}.
\end{align*}

Notice that the last expression  has a similar structure than the right-hand side of (\ref{t132}) where
$B^h_{.*}(t,x):= G(t+h-\cdot,x-*)-G(t-\cdot,x-*)$ is replaced by  
$G(t-.,\xbar-\ast)-G(t-.,x-\ast)$. Hence, we can proceed as in the analysis of the time increments to see that
it suffices to deduce an estimate for
\begin{equation*}
\left\Vert \langle G(t-.,x-\ast) \sigma(u(.,\ast)), \sigma(u(.,\ast))\left[G(t-.,\xbar-\ast)-G(t-.,x-\ast)\right]\rangle_{\mathcal{H}_T}\right\Vert_{3,p},
\end{equation*}
for any $p\in[1,\infty[$.

To pursue the proof, we split the argument of the above expression into two terms
\begin{align*}
S_{3,1}(t,x,\xbar)&=\int_0^t dr \int_{\IR^3} d\xi \int_{\IR^3} d\eta G(t-r,x-\xi)\\
&\quad \times\sigma(u(r,\xi))f(\xi-\eta)\sigma(u(r,\eta))G(t-r,\xbar-\eta)\\
S_{3,2}(t,x,\xbar)&=\int_0^t dr \int_{\IR^3} d\xi \int_{\IR^3} d\eta G(t-r,x-\xi)\\
&\quad\times\sigma(u(r,\xi))f(\xi-\eta)\sigma(u(r,\eta))G(t-r,x-\eta),
\end{align*}
and we apply the change of variables $(\xi\to x-\xi, \eta\to\xbar-\eta)$, $(\xi\to x-\xi, \eta\to x-\eta)$, respectively.  We obtain
\begin{align}
\label{x1}
&S_{3,1}(t,x,\xbar)-S_{3,2}(t,x,\xbar)\nonumber\\
&\quad=\int_0^t dr \int_{\IR^3} \int_{\IR^3}  G(t-r,d\xi ) G(t-r,d\eta)f(\xi-\eta)\sigma(u(r,x-\xi))\nonumber\\
&\quad \quad\times\left(\sigma(u(r,\bar x-\eta))-\sigma(u(r,x-\eta))\right)\nonumber\\
& + \int_0^t dr \int_{\IR^3} \int_{\IR^3}  G(t-r,d\xi ) G(t-r,d\eta)\left[f(x-\bar x-(\xi-\eta))-f(\xi-\eta)\right]\nonumber\\
& \quad \quad\times\sigma(u(r,x-\xi))\sigma(u(r,\bar x-\eta))
 \end{align}
 By Minkowski's inequality the $\Vert\cdot\Vert_{k,p}$--norm of the first term in the right-hand side of (\ref{x1}) is bounded by 
\begin{align}
\label{x2}
%&\left\Vert S_{3,1}(t,x,\xbar)-S_{3,2}(t,x,\xbar)\right\Vert_{3,p}\\
%&\quad \le 
&\int_0^t dr \int_{\IR^3}\int_{\IR^3} G(t-r,d\xi) G(t-r,d\eta) f(\xi-\eta)\nonumber\\
&\qquad \times\left\Vert \sigma(u(r,x-\xi))\left[\sigma(u(r,\xbar-\eta))-\sigma(u(r,x-\eta))\right]\right\Vert_{k,p}\nonumber\\
&\quad\le C |x-\xbar|^{\alpha},
\end{align}
where the very last upper bound follows from (\ref{0.7}), (\ref{hd}) and (\ref{norms}).

For the second term of the right-hand side of (\ref{x1}) we apply Lemma 6.1 of \cite{dss} which implies
\begin{equation*}
\sup_{t\in[0,T]} \int_0^t dr \int_{\IR^3}\int_{\IR^3} G(r,d\xi) G(r,d\eta)\vert f(x-\bar x-(\xi-\eta))-f(\xi-\eta)\vert
\le C |x-\bar x|^{\tilde\alpha},
\end{equation*}
with $\tilde\alpha\in ]0,(2-\beta)\wedge1[$. From this and the properties (\ref{bound1}), (\ref{norms}), we obtain
the upper bound $C |x-\bar x|^{\tilde\alpha}$.

Hence we conclude
\beq
\label{1.16}
\sup_{t\in[0,T]} S_{3}(t,x,\xbar)\le C |x-\xbar|^\alpha,
\eeq
for any $\alpha\in ]0, \gamma_1 \wedge \gamma_2 \wedge \frac{2-\beta}{2}\wedge\frac{1+\delta}{2}[$.
With (\ref{1.14})--(\ref{1.16}), we have proved (\ref{s}).
\smallskip 

\noindent{\bf The rest term in the space expansion}

The contribution of the second order term in (\ref{1.13}) comes from the estimate
\beqn
\left|E\left[g^{\prime\prime}(\hat u(t,x,\xbar))\left(u(t,\xbar)-u(t,x)\right)^2 \right]\right| \le C \Vert g''\Vert_{\infty} h^{2\alpha },
\eeqn
which is a consequence of (\ref{0.7}).

Hence, we have proved that for any fixed $y\in\IR^3$ the mapping $x\in D\mapsto p_{t,x}(y)$ is H\"older continuous of degree
$\alpha\in]0,\gamma_1\wedge\gamma_2\wedge \frac{2-\beta}{2}\wedge\frac{1+\delta}{2}[$, uniformly in
$t\in]0,T]$. 

The proof of the theorem is now complete
 \hfill $\blacksquare$
\medskip

\noindent\textbf{Acknowledgements.}
This paper has been written when the author was visiting the Institute Mittag-Leffler
in Djursholm (Sweden) during a semester devoted to SPDEs.
She would like to express her gratitude for the inspiring environment, the very kind hospitality and the 
financial support provided by this institution.

%%%%%%%%%%%%%%%%%%
%%%%%%%%%%%%%%%%%%
%%%%%%%%%%%%%%%%%%
%%%%%%%%%%%%%%%%%%

%%%%%%%%%%%%%%%%
%%%%%%%%%%%%%%%%
%%%%%%%%%%%%%%%%
%%%%%%%%%%%%%%%%REFERENCES

\end{document}